\documentclass[oneside]{amsart}
\usepackage{float}
\usepackage{amsmath}
\usepackage{graphicx}
\usepackage{amsfonts}
\usepackage{amssymb}
\theoremstyle{plain}
\newtheorem{Theorem}{Theorem}[section]
\newtheorem*{Theorem1}{Theorem 1} 
\newtheorem*{Theorem2}{Theorem 2} 
\newtheorem*{Theorem3}{Theorem 3}
\newtheorem*{Prop23}{Proposition 2.3} 

\newtheorem{Proposition}[Theorem]{Proposition}
\newtheorem{Lemma}[Theorem]{Lemma}

\theoremstyle{remark}

\newtheorem{Remark}{Remark}[section]

\theoremstyle{definition}

\numberwithin{equation}{section}

\begin{document}
\title[Selberg-like Integrals]{On the evaluation of some Selberg-like integrals}
\author{B. Binegar}
\address{Department of Mathematics\\
Oklahoma State University\\
Stillwater, Oklahoma 74078}
\email{binegar@okstate.edu}
\thanks{The author thanks Leticia Barchini for suggesting the problem from which this 
work evolved, and for several effectual suggestions during
the course of this investigation. He also gratefully acknowledges discussions with
Hiroyuki Ochiai and Roger Zierau}
\date{September 28, 2006}
\maketitle
\begin{abstract}Several methods of evaluation are presented for a family of Selberg-like
integrals that arose in the computation of the algebraic-geometric degrees of
a family of multiplicity-free nilpotent $K_{\mathbb{C}}$-orbits. First,
adapting the technique of Nishiyama, Ochiai and Zhu, we present an explicit
evaluation in terms of certain iterated sums over permutations groups.
Secondly, using the theory of symmetric functions we obtain an evaluation as a
product of polynomial of fixed degree times a particular product of gamma
factors (thereby identifying the asymptotics of the integrals with respect to
their parameters). Lastly, we derive a recursive formula for evaluation of 
another general class of Selberg-like integrals, by applying some of the technology 
of generalized hypergeometric functions.
\end{abstract}

\section{Introduction}

In 1944, while still a high school student, Atle Selberg published the
following multivariate generalization of Euler's beta integral formula:
\begin{eqnarray}
&  \int_{\left[  0,1\right]  ^{n}}\left(  \prod_{i=1}^{n}\left(  x_{i}\right)
^{r-1}\left(  1-x_{i}\right)  ^{s-1}\right)  \left(  \prod_{1\leq i<j\leq
n}\left|  x_{i}-x_{j}\right|  ^{2\kappa}\right)  d^{n}x\nonumber\\
& \qquad = \; \prod_{i=1}^{n}\frac{\Gamma\left(  i\kappa+1\right)  \Gamma\left(
r+\left(  i-1\right)  \kappa\right)  \Gamma\left(  s+\left(  i-1\right)
\right)  }{\Gamma\left(  \kappa+1\right)  \Gamma\left(  r+s+2+\left(
n-i-2\right)  \kappa\right)  } \quad .%
\end{eqnarray}
For some 35 years this result laid in deep
hibernation - until it was rediscovered and
vigorously reanimated by Askey \cite{As}, Macdonald \cite{Mac1}, Kor\/anyi
\cite{Kor} and many
others. Indeed, since its rediscovery, the Selberg formula has been
generalized in several directions (\cite{Ao}, \cite{Kad}, \cite{Kan},
\cite{Ri}), and has found important applications in both pure mathematics
(\cite{Kor},\cite{FK}, \cite{KO}, \cite{NO}, \cite{Meh}) and physics
(\cite{Fo},\cite{Ve}).

Still today the Selberg integrals lay at the heart of a fascinating nexus of
representation theory, algebraic geometry, analysis, and combinatorics. A very nice
illustration of this is found in a recent paper by K. Nishiyama and H. Ochiai \cite{NO}.
Most ostensibly, this paper deals with the problem of calculating the Bernstein 
degrees of singular highest weight representations of the metaplectic group. 
However, it turns out that the crux of the matter is to
evaluate an integral of the form
\begin{equation}
\int_{S_{i}}\left(  \prod_{1\leq j\leq i}x_{j}\right)  ^{s}\left(
\prod_{1\leq j<k\leq i}\left(  x_{j}-x_{k}\right)  \right)  ^{d}d^{i}x\quad .
\label{1}%
\end{equation}
This they do by reinterpreting the integral as an integral over the
symmetric cone associated to the space of positive, real symmetric matrices.
They are then able to use results of Faraut and Koranyi \cite{FK} (which go
back to the original Selberg formula via \cite{Kor}) to obtain an explicit evaluation
as a certain product of gamma functions of the parameters.
The authors then remark that their computation of the 
Bernstein degree is equivalent to computing the algebraic-geometric degree of the
determinantal varieties $Sym_{n}\left(  m\right)  =\left\{  X\in M_{n}\left(
\mathbb{C}\right)  \mid\;\right. ~^{t}X=X\ ,\left.
\text{rank}\left(  X\right)  \leq m\right\}  $, and that they have thereby reproduced the
classical formulae of Giambelli.

\bigskip

Inspired by the methods and results of the current generation of Japanese
representation theorists (e.g., \cite{NOTYK}, \cite{NOZ}, \cite{KO}) we have in \cite{B} 
derived a formula for the leading term of the Hilbert polynomials of a family of
multiplicity-free $K_{\mathbb{C}}$-orbits\footnote{We mean $K_{\mathbb{C}}$ orbits
in $\frak{p}$, where $\frak{g} = \frak{k}+ \frak{p}$ is a Cartan decomposition of 
the complexified Lie algebra of $G$.} for any real simple noncompact
group $G$ of classical type. This formula, together with the restricted root data,
reduces the problem of the determining the algebraic-geometric degree of
such an orbit to the evaluation of an integral that is either of the form (\ref{1}),
which happens only when $G/K$ is Hermitian symmetric, or of the form
\begin{equation}
\int_{S_{i}}\left(  \prod_{1\leq j\leq i}x_{j}\right)  ^{s}\left(
\prod_{1\leq j<k\leq i}\left(  x_{j}^{2}-x_{k}^{2}\right)  \right)  ^{d}d^{i}x
\label{2}%
\end{equation}
and where, in both cases, $\mathcal{S}_{i}$ is the domain%
\begin{equation}
\mathcal{S}_{i}=\left\{  x\in\mathbb{R}^{i}\mid x_{1}\geq x_{2}\geq\cdots\geq
x_{i}\geq0\quad,\quad\sum_{j=1}^{i}x_{j}\leq1\right\}  \quad . %
\end{equation}

Let us remark here that integrals of the form (\ref{1}), while replete 
with applications and interpretations (such as in \cite{NO}), can also be 
evaluated simply by a change of variables (\cite{Mac2} pg. 286) and an 
application of the original Selberg formula.   On the other hand, integrals 
of the form (\ref{2}) also have interpretations
and applications outside our particular representation theoretical context. 
For example, integrals of this form arise quite naturally in the 
context of analysis on symmetric domains (see, e.g.,  \cite{FK}, pg. 197). 
Yet it seems, except for very particular values of the parameters,
no explicit evaluation is known.

At first this may seem a modest oversight in the literature; as apparently
the only difference between the two integrals is the replacement of the
linear difference factors linear difference factors $x_{i}-x_{j}$ in (\ref{1})
with the quadratic difference factors quadratic differences $x_{i}^{2}-x_{j}^{2}$.
Morever, from a representation theoretical point of view, the difference/connection 
between the two integrals is just a reflection of the fact that the 
restricted root systems in the Hermitian symmetric case are of type 
$A_{n-1}\approx\left\{  \pm \left( e_{i}- e_{j} \right)\mid1\leq i<j\leq n\right\}  $, 
while, in the non-Hermitian symmetric situation, the
restricted roots systems will involve both sums and differences (as
well as multiples) of the standard basis vectors $e_{i}$. However, the
substitution
$x_{i}-x_{j}\rightarrow x_{i}^{2}-x_{j}^{2},$ also leads to a breaking of a subtle
$x_{i}\longleftrightarrow1-x_{i}$ symmetry in the original Selberg integral. This
loss of symmetry at least partially explains why the integral (\ref{2}) is such a bugger. 

\bigskip

Henceforth, we shall denote by $I_{n,d,p}$ the integral
\begin{equation}
I_{n,d,p} \equiv \int_{\mathcal{S}_{n}}\left(  \prod_{1\leq i\leq n}x_{i}\right)  ^{p}\left(
\prod_{1\leq i<j\leq n}\left(  x_{i}^{2}-x_{j}^{2}\right)  \right)  ^{d}d^{n}x
\label{Indp}%
\end{equation}
where the domain $\mathcal{S}_{n}$ is
\begin{equation}
\mathcal{S}_{n}=\left\{  x\in\mathbb{R}^{n}\mid x_{1}\geq x_{2}\geq\cdots\geq
x_{n}\geq0\quad,\quad\sum_{i=1}^{n}x_{i}\leq1\right\}  \label{6}%
\end{equation}
as above. Note that when $d$ is an even integer, the integrand is a
homogeneous symmetric polynomial and that (in any case) $\mathcal{S}_{n}$ is a
fundamental domain for the action of the symmetric group $\frak{S}_{n}$.
Because of these circumstances we can compute $I_{n,2k,p}$ as $n!$ times the
integral of the same integrand over the simpler region
\[
\Omega_{n}=\frak{S}_{n}\cdot\mathcal{S}_{n}=\left\{  (x_{1},\ldots,x_{n}%
)\in\mathbb{R}^{n}\mid x_{1}\geq0,\;x_{2}\geq0,\;\cdots\;,\;x_{n}\geq
0\;;\sum_{i=1}^{n}x_{i}\leq1\right\}  \quad.
\]
This observation and the method of Nishiyama, Ochiai and Zhu \cite{NOZ} allows
us to obtain in \S 2 the following result.
\begin{Theorem1}
If $d$ is even then
\begin{eqnarray*}
I_{n,d,p}&=&\frac{1}{n!}\frac{1}{\Gamma\left(  a+1\right)  }\sum_{\sigma_{1}%
\in\frak{S}_{n}}\cdots\sum_{\sigma_{d}\in\frak{S}_{n}}\left[sgn\left(  \sigma
_{1}\cdots\sigma_{d}\right) \right. \\
&\times& \left.  \prod_{i=1}^{n}\Gamma\left(  2\sigma_{1}\left(
i\right)  +\cdots+2\sigma_{d}\left(  i\right)  -2d+p+1\right) \right]
\end{eqnarray*}
where 
\[
a \equiv n \left(p +1 \right) + dn \left(n-1 \right)/2 \quad .
\]
\end{Theorem1}

Although this result is fairly explicit, it has two unpleasant features.
First of all, it is valid only for even $d$ (when $d$ is odd the integrand is
skew-symmetric and so the integral can not be computed by extending
its domain to $\Omega_{n}$). Secondly, the sums over permutations are
extremely arduous to compute even for relatively small values for $n$ and $d$.
Toward the end of \S2 we derive the following result which sheds a little more
light on this valuation of $I_{n,d,p}$ for large values of the parameters:

\begin{Prop23} Let $S_{n,d,p}$ be
\[
\sum_{\sigma_{1}%
\in\frak{S}_{n}}\cdots\sum_{\sigma_{d}\in\frak{S}_{n}}\left[sgn\left(  \sigma
_{1}\cdots\sigma_{d}\right) 
 \prod_{i=1}^{n}\Gamma\left(  2\sigma_{1}\left(
i\right)  +\cdots+2\sigma_{d}\left(  i\right)  -2d+p+1\right) \right]
\]
(i.e., $S_{n,d,p}$ is the sum on the right hand side of the formula in Proposition 1.3),
then
\[
S_{n,d,p} =\Phi_{n,d}\left(  p\right)  \prod_{i=1}^{n}%
\Gamma\left(  p+1+d\left(  i-1\right)  \right)
\]
where%
\[
\Phi_{n,d}\left(  p\right)  = \sum_{\sigma\in\frak{S}_{n}}\cdots
\sum_{\rho_{d}\in\frak{S}_{n}}sgn\left(  \sigma_{1}\cdots\sigma_{d}\right)
\prod_{i=1}^{n}\left(  p+1+d\left(  i+1\right)  \right)  _{2\mu_{i}\left(
\mathbf{\sigma}\right)  } \quad .
\]
\end{Prop23}

Here $\left(  n\right)  _{k}=\left(  n\right)  \left(  n-1\right)
\cdots\left(  n-k-1\right)  $ is the usual Pochhammer symbol, and $\mu
_{i}\left(  \mathbf{\sigma}\right)  =\mu_{i}\left(  \sigma_{1},\ldots
,\sigma_{d}\right)  $ is a particular integer-valued function on $\left(
\frak{S}_{n}\right)  ^{d}$ (corresponding roughly to the height of the $i^{th}$
part of a partition constructed from the $\sigma_{1},\ldots,\sigma_{d}$ above
its minimal possible value.) It is readily seem that the degree of the polynomial
factor $\Phi_{n,d}\left(  p\right)  $ is $\leq\frac{d}{2}n\left(  n-1\right)
$. However, explicit computations reveal that
this bound on the degree of $\Phi_{n,d}(p)$ is far from sharp. This is to be remedied
later in \S3 (see Remark 3.1.)

\bigskip

The theory of symmetric functions provides another point of entry into the
topic of generalized Selberg integrals and, in fact, much of the work on
generalized Selberg integrals during the 1980's and 1990's (\cite{Mac1}, \cite{Ao},
\cite{Kad}, \cite{Ri}) was predicated on this point of view. We have already noted
that $I_{n,d,p}$ is an integral of a symmetric or skew-symmetric homogeneous polynomial 
over a fundamental domain for the natural action of the symmetric
group $\frak{S}_{n}$ on $\mathbb{R}^{n}$. In fact, if we denote by
$e_{n}$ the $n^{th}$ elementary symmetric function
\[
e_{n}\left( \mathbf{x} \right) \equiv \prod_{i=1}^{n}x_{i} \quad , 
\]
by $\Delta^{(n)}$ the Vandermonde
determinant%
\[
\Delta^{(n)}\left(  \mathbf{x} \right)  =\det\left[  x_{i}^{n-j}\right]  _{i,j=1,\ldots
,n}=\sum_{1\leq i<j\leq n}\left(  x_{i}-x_{j}\right)
\]
and by $s_{\delta}$ the Schur symmetric function corresponding to the staircase 
partition
$\delta=\left( n-1,\ldots,1,0\right)  $%
\[
s_{\delta}\left( \mathbf{x} \right) \equiv \frac{\det\left[  x_{i}^{n-j+\delta_{j}}\right]  _{i,j=1,\ldots,n}%
}{\det\left[  x_{i}^{n-j}\right]  _{i,j=1,\ldots,n}}=\sum_{1\leq i<j\leq
n}\left(  x_{i}+x_{j}\right) \quad ,
\]
then then we can express the integral (\ref{Indp}) as
\[
I_{n,d,p}=\int_{\mathcal{S}_{n}}\left(  e_{n}\left(  \mathbf{x} \right)  \right)
^{p}\left(  s_{\delta}\left(  \mathbf{x}\right)  \right)  ^{d}\left|  \Delta^{(n)} \left(
\mathbf{x}\right)  \right|  ^{d}d\mathbf{x} \quad .
\]
In \S3 we use the theory of symmetric functions and an integral formula of
Macdonalds, to arrive at the following formula valid for any positive integer
$d$.

\begin{Theorem2}
For $d \in \mathbb{Z}_{\ge0}$,
\[
I_{n,d,p}= \frac{1}{\Gamma\left(n\left(p+d(n-1)\right)+1\right)}\Phi\left(  p\right)  \prod_{i=1}^{n}\Gamma\left( p+1+d\left(
n-i\right)/2  \right)
\]
where $\Phi\left(  p\right)  $ is a polynomial of degree $\leq$ $\frac{d}%
{4}n\left(  n-1\right)  $ if $d$ is even, \ or of degree $\leq\frac{d}%
{4}n\left(  n-1\right)  -\frac{1}{2}\left[  \frac{n}{2}\right]  $ if $d$ is odd.\footnote{Here,
as in common practice,
$\left[  \frac{k}{2}\right]  $ denotes the integer part
of $k/2$.}
\end{Theorem2}

\bigskip

A third approach to the evaluation of Selberg type integrals is 
demonstrated in a remarkable paper by Kaneko \cite{Kan} (which is in turn is based on
ideas of Aomoto \cite{Ao}). Therein, the evaluation of a family of Selberg type
integrals is carried out by identifying certain recursive formulae within the
family, which in turn lead to certain holonomic systems of PDEs satisfied by
the integrals. Kaneko then shows that the unique analytic solutions to these
PDEs are expressible in terms of generalized hypergeometric series and thereby 
obtains explicit evaluations in terms of generalized
hypergeometric functions.

In \S4 we generalize the integrand of the original problem a bit and focus our
efforts on obtaining recursive formulas for another general class of Selberg-like integrals.
More explicitly, we consider integrals of the form
\begin{equation}
J_{n,\kappa}\left( \Phi_{\lambda} \right)  \equiv
\int_{\mathcal{S}_{n}} \Phi_{\lambda} \left(
\mathbf{x}\right) \left(\Delta^{(n)}(\mathbf{x}) \right)^{\kappa}  d\mathbf{x} \label{11}% 
\end{equation}
where $\left\{  \Phi_{\lambda}\right\}  $ is some basis for homogeneous
symmetric polynomials indexed of partitions $\lambda$ (e.g., the
$\Phi_{\lambda}$ might be monomial symmetric functions, or Jack symmetric
functions). We first introduce a change of variables and show how by a sequence of
such transformations any integral of the form (\ref{11}) can, in principle, be 
reduced to a sum of products of beta integrals. 
We then introduce generalized Taylor coefficients $c_{\lambda \mu}$ and
generalized binomial coefficients
$\left(
\begin{array}
[c]{c}%
\lambda\\
\mu
\end{array}
\right)$  
, defined, respectively by the formulas\footnote{
Here $\mathbf{x}_{(n-1)}$ denotes the vector in $\mathbb{R}^{n-1}$ obtained
from $\mathbf{x}_{(n)} = (x_{1},\ldots,x_{n}) \in \mathbb{R}^{n}$ by dropping the last 
coordinate $x_{n}$, and $\mathbf{1}_{(n)} = (1,\ldots,1)\in \mathbb{R}^{n}$.}
\[
\Phi_{\lambda}  \left(  \mathbf{x}_{\left(  n\right)
}\right)  =\sum_{i=0}^{\left|  \lambda\right|  }\sum_{\left|  \mu\right|
=\left|  \lambda\right|  -i}c_{\lambda\mu}\Phi_{\mu}
\left(  \mathbf{x}_{\left(  n-1\right)  }\right)  x_{n}^{i} %
\]
and
\[
\Phi_{\mathbf{\lambda}}\left(  \mathbf{x}_{\left(
n\right)  }+t\mathbf{1}_{\left(  n\right)  }\right)  =\sum_{\mathbf{\mu}%
}\left(
\begin{array}
[c]{c}%
\lambda\\
\mu
\end{array}
\right)  \Phi_{\mathbf{\mu}}\left(  \mathbf{x}_{\left(
n\right)  }\right)  t^{\left|  \lambda\right|  -\left|  \mu\right|  }%
\]
and obtain, from our change of
variables formula, the following recursive formula
\begin{Theorem3}
\begin{eqnarray*}
J_{n,\kappa}\left(  \Phi_{\lambda}\right)
&=&\sum
_{i=0}^{\left|  \lambda\right|  }\sum_{\left|  \mu\right|  =\left|
\lambda\right|  -i}\sum_{\nu}\left(  n\right)  ^{\left|  \nu\right|
-n-i-1}c_{\lambda\mu}\left(
\begin{array}
[c]{c}%
\mu\\
\nu
\end{array}
\right) \\
&&\times  B\left(  \left|  \lambda\right|  +\left|  \nu\right|  +1,n-1+\frac
{\kappa}{2}n\left(  n-1\right)  +\left|  \nu\right|  \right)  J_{n-1,\kappa
}\left(  \Phi_{\nu}\right) \quad .
\end{eqnarray*}
\end{Theorem3}

\bigskip

\section{An explicit evaluation of $I_{n,d,p}$ for the case of even $d$}

In this section we follow the method of Nishiyama, Ochiai and Zhu \cite{NOZ}
to find a closed expression for $I_{n,d,p}$ for case when $d$ is an even integer.

We begin by noting that, for even $d$, the integrand is a symmetric polynomial
in the variables $x_{i}$ and that the region of integration
\[
\mathcal{S}_{n}=\left\{  \mathbf{x}\in\mathbb{R}^{n}\mid x_{1}\geq x_{2}%
\geq\cdots\geq x_{n}\geq0\quad,\quad\sum_{i=1}^{n}x_{i}\leq1\right\}
\]
is a fundamental domain for the natural action of the symmetric group
$\frak{S}_{n}$ on $\mathbb{R}^{n}$. Because of this we can write%
\[
I_{n,d,p}=\frac{1}{n!}\int_{\Omega_{n}}\left(  \prod_{i}x_{i}^{p}\right)
\left(  \prod_{1\leq i<j\leq n}\left(  x_{i}^{2}-x_{j}^{2}\right)
^{d}\right)  d^{n}x
\]
where $\Omega_{n}$ is the much simpler simplex
\[
\Omega_{n}=\frak{S}_{n}\cdot\mathcal{S}_{n}=\left\{  \mathbf{x}\in
\mathbb{R}^{n}\mid x_{i}\geq0\ ,\ i=1,\ldots n\ ,\quad\text{and}\quad
\sum_{i=1}^{n}x_{i}\leq1\right\} \quad .
\]
Next, we set%
\[
a=n\left(  p+d\left(  n-1\right)  +1\right)
\]
and note that degree of the integrand is $a-n$.

We now use the identity%
\[
\Gamma\left(  a+1\right)  \int_{\Omega_{n}}f\left(  \mathbf{x}\right)
d\mathbf{x}=\int_{\Omega_{n}\times\left(  0,\infty\right)  }f\left(
\mathbf{x}\right)  s^{a}e^{-s}d\mathbf{x}ds
\]
and make a change of variables%
\[%
\begin{array}
[c]{l}%
y_{i}=sx_{i}\\
~\\
t=s\left(  1-\sum_{i=1}^{n}x_{i}\right)
\end{array}
\quad\Longleftrightarrow\quad%
\begin{array}
[c]{l}%
x_{i}=\frac{y_{i}}{t+\sum_{k=1}^{n}y_{k}}\\
~\\
s=t+\sum_{k=1}^{n}y_{k}%
\end{array}
\]
which maps $\Omega_{n}\times\left(  0,\infty\right)  $ diffeomorphically onto
$\left(  0,\infty\right)  ^{n}\times\left(  0,\infty\right)  $. The Jacobian
of this transformation is easily seen to be
\[
\frac{\partial\left(  y,t\right)  }{\partial\left(  x,s\right)  }%
=s^{n}=\left(  t+\sum_{k=1}^{n}y_{k}\right)  ^{n}%
\]
and so we have, for any function $f$ homogeneous of degree $a-n$,%
\begin{eqnarray*}
\Gamma\left(  a+1\right)  \int_{\Omega_{n}}f\left(  \mathbf{x}\right)
d\mathbf{x}  &  =&\int_{\Omega_{n}\times\left(  0,\infty\right)  }f\left(
\mathbf{x}\right)  s^{a}e^{-s}d\mathbf{x}ds\\
&  =&\int_{\left(  0,\infty\right)  ^{n}\times\left(  0,\infty\right)
}f\left(  \frac{\mathbf{y}}{t+\sum_{k=1}^{n}y_{k}}\right)  \left(
t+\sum_{k=1}^{n}y_{k}\right)  ^{a} \\
&&\times \quad e^{-t-\sum_{k=1}^{n}y_{k}}\left(
t+\sum_{k=1}^{n}y_{k}\right)  ^{n}d\mathbf{y}dt\\
&  =&\int_{\left(  0,\infty\right)  ^{n}\times\left(  0,\infty\right)
}f\left(  \mathbf{y}\right)  e^{-t-\sum_{k=1}^{n}y_{k}}d\mathbf{y}dt\\
&  =&\int_{\left(  0,\infty\right)  ^{n}}f\left(  \mathbf{y}\right)
e^{-\sum_{k=1}^{n}y_{k}}d\mathbf{y} \quad . %
\end{eqnarray*}

We thus arrive at%
\[
I_{n,d,p}=\frac{1}{n!}\frac{1}{\Gamma\left(  a+1\right)  }\int_{\left(
0,\infty\right)  ^{n}}\left(  \prod_{i=1}^{n}y_{i}^{p}\right)  \left(
\prod_{1\leq i<j\leq n}\left(  y_{i}^{2}-y_{j}^{2}\right)  ^{d}\right)
e^{-\sum_{k=1}^{n}y_{k}}d\mathbf{y} \quad .%
\]

The next step is to expand the second product using the identity%
\[
\prod_{1\leq i<j\leq n}\left(  y_{i}^{2}-y_{j}^{2}\right)  =\det\ \left(
y_{j}^{2\left(  i-1\right)  }\right)  _{\substack{i=1,\ldots,n\\j=1,\ldots
,n}}=\sum_{\sigma\in\frak{S}_{n}}sgn\left(  \sigma\right)  \prod_{i=1}%
^{n}y_{i}^{2\left(  \sigma\left(  i\right)  -1\right)  } \quad . %
\]
We have%
\begin{eqnarray*}
\left(  \sum_{\sigma\in\mathbb{S}_{n}}sgn\left(  \sigma\right)  \prod
_{i=1}^{n}y_{i}^{2\left(  \sigma\left(  i\right)  -1\right)  }\right)  ^{d}
&  =&\sum_{\sigma_{1}\in\frak{S}_{n}}\cdots\sum_{\sigma_{d}\in\frak{S}_{n}%
}\left[ sgn\left(  \sigma_{1}\right)  \cdots sgn\left(  \sigma_{d}\right)\right.\\
&& \left. \quad \times 
\prod_{i=1}^{n}y_{i}^{2\left(  \sigma_{1}\left(  i\right)  -1\right)
+\cdots+2\left(  \sigma_{d}\left(  i\right)  -1\right)  } \right]\\
&  =&\sum_{\sigma\in\frak{S}_{n}}\cdots\sum_{\sigma_{dd}\in\frak{S}_{n}%
}sgn\left(  \sigma_{1}\cdots\sigma_{d}\right)  \prod_{i=1}^{n}y_{i}^{2\left(
\sum_{j=1}^{d}\sigma_{j}\left(  i\right)  \right)  -2d}%
\end{eqnarray*}
and so%
\begin{eqnarray*}
I_{n,d,p}  &=&\frac{1}{n!}\frac{1}{\Gamma\left(  a+1\right)  }\sum
_{\sigma_{1}\in\frak{S}_{n}}\cdots\sum_{\sigma_{d}\in\frak{S}_{n}}sgn\left(
\sigma_{1}\cdots\sigma_{d}\right)  \int_{\left(  0,\infty\right)  ^{n}} \\
&&\times \prod_{i=1}^{n}y_{i}^{2\left(  \sum_{j=1}^{d}\sigma_{j}\left(  i\right)
\right)  -2d+p}e^{-\sum_{k=1}^{n}y_{k}}d\mathbf{y}\\
&=&\frac{1}{n!}\frac{1}{\Gamma\left(  a+1\right)  }\sum_{\sigma_{1}%
\in\frak{S}_{n}}\cdots\sum_{\sigma_{d}\in\frak{S}_{n}}\left[ sgn\left(  \sigma
_{1}\cdots\sigma_{d}\right) \right. \\
&& \left. \prod_{i=1}^{n}\left(  \int_{0}^{\infty}%
y_{i}^{2\left(  \sum_{j=1}^{d}\sigma_{j}\left(  i\right)  \right)
-2d+p}e^{-y_{i}}dy_{i}\right) \right] \\
& =&\frac{1}{n!}\frac{1}{\Gamma\left(  a+1\right)  }\sum_{\sigma_{1}%
\in\frak{S}_{n}}\cdots\sum_{\sigma_{d}\in\frak{S}_{n}}\left[ sgn\left(  \sigma
_{1}\cdots\sigma_{d}\right) \right. \\
&& \quad \times
\left. \prod_{i=1}^{n}\Gamma\left(  2\sigma_{1}\left(
i\right)  +\cdots+2\sigma_{d}\left(  i\right)  -2d+p+1\right) \right]
\end{eqnarray*}
where we have used the Euler's formula for the gamma function%
\[
\Gamma\left(  x\right)  =\int_{0}^{\infty}t^{x-1}e^{-t}dt\quad.
\]

In summary,

\begin{Theorem1}
If $d$ is even and
\[
I_{n,d,p}\equiv\int_{\mathcal{S}_{n}}\left(  \prod_{i}x_{i}^{p}\right)
\left(  \prod_{1\leq i<j\leq n}\left(  x_{i}^{2}-x_{j}^{2}\right)
^{d}\right)  d^{n}x
\]
then%
\begin{eqnarray*}
I_{n,d,p}&=&\frac{1}{n!}\frac{1}{\Gamma\left(  a+1\right)  }\sum_{\sigma_{1}%
\in\frak{S}_{n}}\cdots\sum_{\sigma_{d}\in\frak{S}_{n}}\left[ sgn\left(  \sigma
_{1}\cdots\sigma_{d}\right) \right.\\
&& \times \left. \prod_{i=1}^{n}\Gamma\left(  2\sigma_{1}\left(
i\right)  +\cdots+2\sigma_{d}\left(  i\right)  -2d+p+1\right) \right]
\end{eqnarray*}
where $a=n\left(  p+d\left(  n-1\right)  +1\right)  $.
\end{Theorem1}

\bigskip

We now focus our attention on the product of gamma factors%
\[
\prod_{i=1}^{n}\Gamma\left(  2\sigma_{1}\left(  i\right)  +\cdots+2\sigma
_{d}\left(  i\right)  -2d+p+1\right) \quad . 
\]
We begin by forming the vector%
\[
\left[  \sum_{j=1}^{d}\sigma_{j}\left(  1\right)  ,\sum_{j=1}^{d}\sigma
_{j}\left(  2\right)  ,\ldots,\sum_{j=1}^{d}\sigma_{j}\left(  n\right)
\right]
\]
and then reordering the components\ in\ increasing order to form
\[
\mathbf{\gamma}\left(  \mathbf{\sigma}\right)  \equiv\left[  \min_{i}\left\{
\sum_{j=1}^{d}\sigma_{j}\left(  i\right)  \right\}  ,\ldots,\max_{i}\left\{
\sum_{j=1}^{d}\sigma_{j}\left(  i\right)  \right\}  \right] \quad .
\]

\begin{Lemma}
For any arrangement $\mathbf{\sigma}\in\left(  \frak{S}_{n}\right)  ^{d}$ we
have%
\[
\gamma_{i}\left(  \mathbf{\sigma}\right)  \geq\frac{d}{2}\left(  i+1\right) \quad .
\]
\end{Lemma}

\emph{Proof.} We first note that%
\[
\gamma_{1}\left(  \mathbf{\sigma}\right)  \geq d
\]
follows readily from the requirement that each $\sigma_{j}\left(  i\right)
\geq1$. \ Next, we note that for any arrangement $\left(  \sigma_{1}%
,\ldots,\sigma_{d}\right)  $%
\[
\sum_{j=1}^{d}\sum_{i=1}^{n}\sigma_{i}\left(  j\right)  =\sum_{j=1}^{d}\left(
1+2+\cdots+n\right)  =\frac{d}{2}n\left(  n+1\right)
\]
and that the particular arrangement where
\[
\sigma_{k}=\left\{
\begin{array}
[c]{ll}%
\left[  1,2,\ldots,n-1,n\right]  & \text{if }k\text{ is even}\\
\left[  n,n-1,\ldots,2,1\right]  & \text{if }k\text{ is odd}%
\end{array}
\right.
\]
leads to
\begin{align*}
\mathbf{\gamma}\left(  \mathbf{\sigma}\right)   &  =\left[  \frac{d}{2}\left(
1\right)  +\frac{d}{2}\left(  n\right)  ,\frac{d}{2}\left(  2\right)
+\frac{d}{2}\left(  n-1\right)  ,\ldots,\frac{d}{2}\left(  n\right)  +\frac
{d}{2}\left(  1\right)  \right] \\
&  =\left[  \frac{d}{2}\left(  n+1\right)  ,\frac{d}{2}\left(  n+1\right)
,\ldots,\frac{d}{2}\left(  n+1\right)  \right] \quad .
\end{align*}
Now note that one cannot decrease the last component of $\mathbf{\gamma
}\left(  \mathbf{\sigma}\right)  $ further without violating the requirement
that $\sum_{i=1}^{n}\gamma_{i}\left(  \mathbf{\sigma}\right)  =\frac{d}%
{2}n\left(  n-1\right)  $ (and the stipulated ordering $\gamma_{i}\left(
\mathbf{\sigma}\right)  \leq\gamma_{i+1}\left(  \mathbf{\sigma}\right)  $).
Thus, we have%
\[
\gamma_{n}\left(  \mathbf{\sigma}\right)  \geq\frac{d}{2}\left(  n+1\right) \quad .
\]
Finally, we observe that the arrangement%
\[
\sigma_{k}=\left\{
\begin{array}
[c]{ll}%
\left[  1,2,\ldots,i-1,i,i+1,\ldots,n\right]  & \text{if }k\text{ is even}\\
\left[  i,i-1,\ldots,2,1,i+1,\ldots,n\right]  & \text{if }k\text{ is odd}%
\end{array}
\right.
\]
leads to%
\[
\mathbf{\gamma}\left(  \mathbf{\sigma}\right)  \mathbf{=}\left[  \frac{d}%
{2}\left(  i+1\right)  ,\ldots,\frac{d}{2}\left(  i+1\right)  ,d\left(
i+1\right)  ,\ldots,d\left(  n\right)  \right]
\]
and so, by essentially the reasoning as above,%
\[
\gamma_{i}\left(  \mathbf{\sigma}\right)  \geq\frac{d}{2}\left(  i+1\right) \quad .
\]
\qed

\begin{Proposition}
Let%
\[
S_{n,d}\left(  p\right)  =\sum_{\sigma\in\frak{S}_{n}}\cdots\sum_{\rho_{d}%
\in\frak{S}_{n}}sgn\left(  \sigma_{1}\cdots\sigma_{d}\right)  \prod_{i=1}%
^{n}\Gamma\left(  2\sigma_{1}\left(  i\right)  +\cdots+2\sigma_{d}\left(
i\right)  -2d+p+1\right) \quad .
\]
Then
\[
S_{n,d}\left(  p\right)  =\Phi_{n,d}\left(  p\right)  \prod_{i=1}^{n}%
\Gamma\left(  p+1+d\left(  i-1\right)  \right)
\]
where%
\begin{equation}
\Phi_{n,d}\left(  p\right)  =\left(  \sum_{\sigma\in\frak{S}_{n}}\cdots
\sum_{\rho_{d}\in\frak{S}_{n}}sgn\left(  \sigma_{1}\cdots\sigma_{d}\right)
\prod_{i=1}^{n}\left(  p+1+d\left(  i+1\right)  \right)  _{2\mu_{i}\left(
\sigma\right)  }\right) \label{2.1}
\end{equation}
is a polynomial in $p$ of degree $\ \leq\frac{d}{2}n\left(  n-1\right)  $.
\end{Proposition}

\emph{Proof.} Set%
\[
\mu_{i}\left(  \mathbf{\sigma}\right)  =\gamma_{i}\left(  \mathbf{\sigma
}\right)  -\frac{d}{2}\left(  i+1\right)
\]
so that%
\[
\mu_{i}\left(  \mathbf{\sigma}\right)  \geq0\qquad,\qquad i=1,\ldots,n \quad .
\]
For each term of the iterated sum we can arrange the gamma factors so that
their arguments are non-decreasing. In other words we can write%
\begin{align*}
S_{n,d}\left(  p\right)   &  =\sum_{\sigma\in\frak{S}_{n}}\cdots\sum_{\rho
_{d}\in\frak{S}_{n}}sgn\left(  \sigma_{1}\cdots\sigma_{d}\right)  \prod
_{i=1}^{n}\Gamma\left(  2\gamma_{i}\left(  \mathbf{\sigma}\right)
-2d+p+1\right) \\
&  =\sum_{\sigma\in\frak{S}_{n}}\cdots\sum_{\rho_{d}\in\frak{S}_{n}}sgn\left(
\sigma_{1}\cdots\sigma_{d}\right)  \prod_{i=1}^{n}\Gamma\left(  2\mu
_{i}\left(  \mathbf{\sigma}\right)  +d\left(  i+1\right)  -2d+p+1\right) \\
&  =\sum_{\sigma\in\frak{S}_{n}}\cdots\sum_{\rho_{d}\in\frak{S}_{n}}sgn\left(
\sigma_{1}\cdots\sigma_{d}\right)  \prod_{i=1}^{n}\Gamma\left(  p+1+d\left(
i-1\right)  +2\mu_{i}\left(  \mathbf{\sigma}\right)  \right) \quad .
\end{align*}
We now introduce the Pochhammer symbols $\left(  k\right)  _{n}$ defined by%
\[
\left(  k\right)  _{n}\equiv\frac{\Gamma\left(  k+n\right)  }{\Gamma\left(
k\right)  } \qquad , %
\]
noting that $\left(  k\right)  _{0}=1$ and that for positive integers $n$%
\[
\left(  k\right)  _{n}=\left(  k\right)  \left(  k+1\right)  \cdots\left(
k+n-1\right) \quad .
\]

Returning to our expression for $S_{n,d}\left(  p\right)  $ we have%
\begin{eqnarray*}
S_{n,d}\left(  p\right)   &  =&\sum_{\sigma\in\frak{S}_{n}}\cdots\sum_{\rho
_{d}\in\frak{S}_{n}}sgn\left(  \sigma_{1}\cdots\sigma_{d}\right)  \prod
_{i=1}^{n}\Gamma\left(  p+1+2\mu_{i}\left(  \sigma\right)  +d\left(
i-1\right)  \right) \\
&  =&\sum_{\sigma\in\frak{S}_{n}}\cdots\sum_{\rho_{d}\in\frak{S}_{n}}sgn\left(
\sigma_{1}\cdots\sigma_{d}\right)  \prod_{i=1}^{n}\left(  p+1+d\left(
i-1\right)  \right)  _{2\mu_{i}\left(  \sigma\right)  }\Gamma\left(
p+1+d\left(  i-1\right)  \right) \\
&=&\left(  \sum_{\sigma\in\frak{S}_{n}}\cdots\sum_{\rho_{d}\in\frak{S}_{n}%
} sgn\left(  \sigma_{1}\cdots\sigma_{d}\right) \prod_{i=1}^{n}\left(
p+1+d\left(  i-1\right)  \right)  _{2\mu_{i}\left(  \sigma\right)  }\right) \\
&&\times \quad \prod_{i=1}^{n}\Gamma\left(  p+1+d\left(  i+1\right)  \right)  \quad .
\end{eqnarray*}
Finally we note that each factor
\[
\prod_{i=1}^{n}\left(  d\left(  i-1\right)  +p+1\right)  _{2\mu_{i}\left(
\sigma\right)  }%
\]
is a polynomial in $p$ of total degree%
\begin{align*}
\sum_{i=1}^{n}2\mu_{i}\left(  \sigma\right)   &  =\sum_{i=1}^{n}\left(
2\gamma_{i}\left(  \mathbf{\sigma}\right)  -d\left(  i+1\right)  \right) \\
&  =2\sum_{i=1}^{n}\gamma_{i}\left(  \mathbf{\sigma}\right)  -d\left(
\frac{1}{2}n\left(  n+1\right)  +n\right) \\
&  =2\sum_{i=1}^{n}\sum_{j=1}^{d}\sigma_{j}\left(  i\right)  -d\left(
\frac{1}{2}n\left(  n+1\right)  +n\right) \\
&  =2d\left(  \sum_{i=1}^{n}i\right)  -d\left(  \frac{1}{2}n\left(
n+1\right)  +n\right) \\
&  =2d\frac{1}{2}n\left(  n+1\right)  -\frac{d}{2}n\left(  n+1\right)  -dn\\
&  =\frac{d}{2}n\left(  n-1\right) \quad .
\end{align*}
We conclude that%
\[
S_{n,d}\left(  p\right)  =\Phi_{n,d}\left(  p\right)  \prod_{i=1}^{n}%
\Gamma\left(  p+1+d\left(  i-1\right)  \right)
\]
with%
\begin{equation}
\Phi_{n,d}\left(  p\right)  =\left(  \sum_{\sigma\in\frak{S}_{n}}\cdots
\sum_{\rho_{d}\in\frak{S}_{n}}sgn\left(  \sigma_{1}\cdots\sigma_{d}\right)
\prod_{i=1}^{n}\left(  p+1+d\left(  i+1\right)  \right)  _{2\mu_{i}\left(
\sigma\right)  }\right) \quad .
\end{equation}
a polynomial of degree $\leq$ $\frac{d}{2}n\left(  n-1\right)  $. \qed

We remark here that the bound $\deg\left(  \Phi_{n,d}\left(  p\right)
\right)  \leq\frac{d}{2}n\left(  n-1\right)  $ is certainly not
optimal. To see this, note that the Pochhammer products%
\[
\prod_{i=1}^{n}\left(  p+1+d\left(  i+1\right)  \right)  _{2\mu_{i}\left(
\sigma\right)  }%
\]
are all monic polynomials, and consequently when we sum over the arrangements in
$\left(  \frak{S}_{n}\right)  ^{d}$ , the $sgn\left(  \sigma_{1}\cdots
\sigma_{d}\right)  $ factors will lead to a complete cancellation of the
leading terms. In fact, explicit computations of the right hand side of
(\ref{2.1}) reveal that at least for small $n$ and $d$ the actual degree of
$\Phi_{n,d}\left(  p\right)  $ is $\frac{d}{4}n\left(  n-1\right)  $; that is,
that that the terms of degree $\frac{d}{2}n\left(  n-1\right)  ,\frac{d}%
{2}n\left(  n-1\right)  -1,\ldots,\frac{d}{4}n\left(  n-1\right)  +1$ all,
quite remarkably, cancel. Unfortunately, we have yet to find a direct
combinatorial argument as to why the first $\frac{d}{4}n\left(  n-1\right)  $
leading terms all cancel. However, in \S3 we shall succeed not only in
extending our results to the case of odd $d$, but also in obtaining a least
upper bound on the degree of the polynomial factor $\Phi_{n,d}\left(
p\right)  $ for arbitrary positive integers $d$ and $n$.

\subsection{The case when $d=2$}

When $d=2$ we can obtain a mor succinct determinantal formula for the polynomial factor
$\Phi_{n,d}\left(  p\right)  $. Starting with%
\[
S_{n,2}\left(  p\right)  =\sum_{\sigma\in\frak{S}_{n}}\sum_{\rho\in
\frak{S}_{n}}sgn\left(  \sigma\right)  sgn\left(  \rho\right)  \prod_{i=1}%
^{n}\Gamma\left(  2\sigma\left(  i\right)  +2\rho\left(  i\right)
-3+p\right)
\]
we can write%
\[
\prod_{i=1}^{n}\Gamma\left(  2\sigma\left(  i\right)  +2\rho\left(  i\right)
+p-3\right)  =\prod_{i=1}^{n}\Gamma\left(  2i+2\rho\left(  \sigma
^{-1}\left(  i\right)  \right)  -3+p+1\right)
\]
and then replace the sum over $\rho\in\frak{S}_{n}$ with the sum over
$\tau=\rho\sigma^{-1}\in\frak{S}_{n}$. \ Noting that
\[
sgn\left(  \tau\right)  =sgn\left(  \rho\right)  sgn\left(  \sigma
^{-1}\right)  =sgn\left(  \rho\right)  sgn\left(  \sigma\right) \quad ,
\]
we obtain in this way%
\begin{align*}
S_{n,2}\left(  p\right)   &  =\sum_{\sigma\in\frak{S}_{n}}\sum_{\tau
\in\frak{S}_{n}}sgn\left(  \tau\right)  \prod_{i=1}^{n}\Gamma\left(
2i+2\tau\left(  i\right)  -3+p\right) \\
&  =n!\sum_{\tau\in\frak{S}_{n}}sgn\left(  \tau\right)  \prod_{i=1}^{n}%
\Gamma\left(  p+2i-1+2\tau\left(  i\right)  -2\right) \\
&  =n!\sum_{\tau\in\frak{S}_{n}}sgn\left(  \tau\right)  \prod_{i=1}^{n}\left(
p-1+2i\right)  _{2\tau\left(  i\right)  -2}\Gamma\left(  p-1+2i\right) \\
&  =n!\det\left[  \left(  p-1+2i\right)  _{2j-2}\right]
_{\substack{i=1,\ldots,n\\j=1,\ldots,n}}\prod_{i=1}^{n}\Gamma\left(
p-1+2i\right)
\end{align*}
and so%
\[
\Phi_{n,2}\left(  p\right)  =n!\det\left[  \left(  p-1+2i\right)
_{2j-2}\right]  _{\substack{i=1,\ldots,n\\j=1,\ldots,n}}\quad.
\]

\begin{Remark}
This formula for $\Phi_{n,2}\left(  p\right)  $ appears as Remark (3.7)\ in [NOZ].
\end{Remark}

\begin{Remark}
The above determinant formula of $\Phi_{n,2}\left(  p\right)  $ seems to
predict a polynomial of total degree%
\[
\sum_{j=1}^{n}\left(  2j-2\right)  =n\left(  n+1\right)  -2n=n\left(
n-1\right)  =\frac{d}{2}n\left(  n-1\right)
\]
in agreement with the upper bound stated in Prop. 2.2.
But, as remarked above, explicit computations reveal that, somewhat miraculously, 
the first
$\left( deg\;\Phi_{n,d} \right)/2$ leading terms all cancel and one has%
\[
\deg\left(  \Phi_{n,2}\left(  p\right)  \right)  =\frac{1}{2}n\left(
n-1\right)  =\frac{d}{4}n\left(  n-1\right) \quad .
\]
\end{Remark}

\section{The Integral for the case $d\in\mathbb{Z}_{>0}$}

\subsection{Symmetric polynomials and an integral formula of Macdonald}

In order to keep the exposition of our results self-contained, we begin 
with a rapid review of the pertinent theory of symmetric functions. A
standard reference for this material is of course \cite{Mac2}.

By a symmetric polynomial of $n$ variables we mean a polynomial $p\in
\mathbb{C}\left[  x_{1},\ldots,x_{n}\right]  $ invariant under the natural of
the symmetric group $\frak{S}_{n}$:
\[
p\left(  x_{1},\ldots,x_{n}\right)  =\left(  \sigma\cdot p\right)  \left(
x_{1},\ldots,x_{n}\right)  \equiv p\left(  x_{\sigma\left(  1\right)  }%
,\ldots,x_{\sigma\left(  n\right)  }\right)  \qquad\forall\ \sigma\in
\frak{S}_{n}\quad.
\]
The action of $\frak{S}_{n}$ preserves the subspaces of $\mathbb{C}\left[
x_{1},\ldots,x_{n}\right]  $ consisting of homogeneous polynomials of fixed
total degree. We shall denote by $\Lambda_{\left(  n\right)  }^{m}$ the space
of homogeneous symmetric polynomials in $n$ variables of total degree $m$, so
that%
\[
\Lambda_{\left(  n\right)  }\equiv\mathbb{C}\left[  x_{1},\ldots,x_{n}\right]
^{\frak{S}_{n}}=\bigoplus_{m=0}^{\infty}\Lambda_{\left(  n\right)  }^{m} \quad .%
\]
There are several fundamental bases for the subspaces $\Lambda_{\left(
n\right)  }^{m}$, each parameterized by partitions of length $n$ and weight
$m$. A partition $\lambda$ of length $n$ is simply a non-increasing list of
non-negative integers; i.e., $\lambda=\left[  \lambda_{1},\lambda_{2}%
,\ldots,\lambda_{n}\right]  $ with $\lambda_{1}\geq\lambda_{2}\geq\cdots
\geq\lambda_{n}\geq0$. The weight $\left|  \lambda\right|  $ of a partition $\lambda$ is
the sum of its parts; i.e., $\left|  \lambda\right|  =\sum_{i=1}^{n}%
\lambda_{i}$. There is a partial ordering of the set of partitions of weight
$w$, called the \emph{dominance partial ordering}, defined as follows
\[
\mu\leq\lambda\quad\Longrightarrow\quad \left| \lambda \right| = \left| \mu \right|
\quad \text{and} \quad \sum_{j=1}^{i}\left(  \lambda_{j}%
-\mu_{j}\right)  \geq0\quad,\quad\text{for }i=1,2,\ldots,n \; .
\]

Let $\lambda=\left[  \lambda_{1},\ldots,\lambda_{n}\right]  $ be a partition,
and let $\mathbf{x}^{\lambda}=x_{1}^{\lambda_{1}}\cdots x_{n}^{\lambda_{n}}$ be the
corresponding monomial. The \emph{monomial symmetric function m}$_{\lambda}$
is the sum of all distinct monic monomials that can be obtained from
$\mathbf{x}^{\lambda}$ by permuting the $x_{i}$'s. Every homogeneous symmetric
polynomial of degree $m$ can be uniquely expressed as a linear combination of
the $m_{\lambda}$ with $\left|  \lambda\right|  =m$.

The \emph{power sum} symmetric polynomials are defined as follows. For each
$r\geq1$, let
\[
p_{r}(\mathbf{x}=m_{\left(  r\right)  }(\mathbf{x} )=\sum_{i=1}^{n}x_{i}^{r}%
\]
and then for any partition $\lambda,$ set
\[
p_{\lambda}\left(  \mathbf{x}\right)  =p_{\lambda_{1}}\left(  x\right)  p_{\lambda_{2}%
}\left(  x\right)  \cdots p_{\lambda_{n}}\left(  x\right) \quad .
\]
The power sum symmetric functions $p_{\lambda}\left(  x\right)  $ with
$\left|  \lambda\right|  =m$ provide another basis for the homogeneous
symmetric polynomials of degree $m$. In what follows, the power sum symmetric
functions are only used to define a particular inner product for the symmetric
polynomials; i.e., an inner product will be defined by specifying matrix
entries with respect to the basis of power sum symmetric polynomials.

The \emph{Schur polynomials} provide yet another basis. These can be defined
as follows. For any partition $\lambda$ of length $n$,
\[
a_{\lambda}\left(  \mathbf{x}\right)  =\det\left(  x_{i}^{\lambda_{j}}\right) \quad .
\]
The $a_{\lambda}\left( \mathbf{ x} \right)  $ are obviously odd with respect to the
action of the permutation group $\frak{S}_{n}$; i.e. $a_{\lambda}\left(
\sigma\left( \mathbf{ x} \right)  \right)  =sgn\left(  \sigma\right)  a_{\lambda}\left(
x\right)  $ for all $\sigma\in\frak{S}_{n}$. However, it turns out that
$a_{\lambda}=0$ for and $\lambda<\delta\equiv\left[  n-1,n-2,\ldots
,1,0\right]  $, and that, for every partition $\lambda$, $a_{\delta}\left(
\mathbf{x} \right)  $ divides $a_{\lambda+\delta}\left(  \mathbf{x}\right)  $. Indeed,
\[
s_{\lambda}\left(  \mathbf{x}\right)  \equiv\frac{a_{\lambda+\delta}\left(  \mathbf{x}\right)
}{a_{\delta}\left(  \mathbf{x}\right)  }%
\]
is a symmetric polynomial of degree $\left|  \lambda\right|  $. The
polynomials $s_{\lambda}\left( \mathbf{ x} \right)  $ are the Schur symmetric
polynomials. The Schur polynomials $\left\{  s_{\lambda}\mid\left|
\lambda\right|  =m\right\}  $ provide another fundamental basis for the
symmetric polynomials that are homogeneous of degree $m$. A special case that
will be important to us later on is
\[
s_{\delta}\left(  \mathbf{x}\right)  =\prod_{1\leq i<j\leq n}\left(  x_{i}%
+x_{j}\right)
\]
where again $\delta=\left[  n-1,n-2,\ldots,1,0\right]  $. We note also the
fact that $a_{\delta}\left( \mathbf{ x} \right)  $ is just the Vandermonde determinant
\[
a_{\delta}\left(  \mathbf{x} \right)  =\det\left[
\begin{array}
[c]{ccccc}%
x_{1}^{n-1} & x_{1}^{n-2} & \cdots &  x_{1} & 1\\
x_{2}^{n-1} & x_{2}^{n-2} & \cdots &  x_{2} & 1\\
\vdots & \vdots & \ddots & \vdots & \vdots\\
x_{n-1}^{n-1} & x_{n-1}^{n-2} & \cdots &  x_{n-1} & 1\\
x_{n}^{n-1} & x_{n}^{n-2} & \cdots &  x_{n} & 1
\end{array}
\right]  =\prod_{1\leq i<j\leq n}\left(  x_{i}-x_{j}\right) = \Delta(\mathbf{x}) \quad.
\]

\emph{Jack's symmetric functions} $P_{\lambda}^{\left(  \alpha\right)  }$ are
symmetric functions indexed by partitions and depending rationally on a
parameter $\alpha$ which interpolate between the Schur functions $s_{\lambda
}\left( \mathbf{ x} \right)  $ , the monomial symmetric functions and two other bases
associated with spherical symmetric polynomials on symmetric spaces.  
They are (uniquely) characterized by two properties

\begin{itemize}
\item [(i)]$P_{\lambda}^{\left(  \alpha\right)  }\left( \mathbf{ x} \right)
=m_{\lambda}\left( \mathbf{ x} \right)  +\sum_{\substack{\mu<\lambda\\\left|
\mu\right|  =\left|  \lambda\right|  }}c_{\lambda,\mu}m_{\mu}\left( \mathbf{ x} \right)
\quad$; that is, the ``leading term'' of $P_{\lambda}^{\left(  \alpha\right)
}$ is the monomial symmetric functions $m_{\lambda}$ and the remaining terms
involve only monomial symmetric functions $m_{\mu}$ for which the partition
index $\mu$ is less than $\lambda$ with respect to the dominance ordering.

\item[(ii)] When one defines a scalar product on the vector space of
homogeneous polynomials of degree $m$ by
\[
\left\langle p_{\lambda},p_{\mu}\right\rangle =\delta_{\lambda,\mu}%
\alpha^{\ell\left(  \lambda\right)  }\prod_{r=1}^{n}\left(  r^{m_{\lambda
}\left(  r\right)  }\cdot m_{\lambda}\left(  r\right)  !\right)
\]
where $m_{\lambda}\left(  r\right)  $ is the number of times the integer $r$
appears in $\lambda$, then
\[
\left\langle P_{\lambda}^{\left(  a\right)  },P_{\mu}^{\left(  \alpha\right)
}\right\rangle =0\qquad\text{if }\mu\neq\lambda\quad.
\]
\end{itemize}
We note that this definition basically ensures that the $P_{\lambda}^{\left(
\alpha\right)  }\left( \mathbf{ x} \right)  $ are constructible via a Gram-Schmidt
process (although not quite straightforwardly, as the dominance ordering is
only a partial ordering).

In fact, the Jack symmetric polynomials $P_{\lambda}^{\left(  \alpha\right)
}$ for special values $\alpha=2,\ 1$, and $\frac{1}{2}$ can also be
characterized as the spherical polynomials for, respectively, $GL\left(
n,\mathbb{R}\right)  /O\left(  n\right)  $, $GL\left(  n,\mathbb{H}\right)
/U\left(  \mathbb{H}\right)  $. When $\alpha=1$ the Jack symmetric polynomials
coincide with the spherical polynomials for $GL\left(  n,\mathbb{C}\right)
/U\left(  n\right)  $, as well the Schur polynomials.

Next we recall the following well known property of Schur polynomials%
\[
s_{\mu}\left( \mathbf{ x} \right)  s_{\nu}\left( \mathbf{ x} \right)  =
\sum_{\lambda\leq\mu+\nu}K_{\mu\nu}^{\lambda}s_{\lambda}%
\]
where the coefficients $K_{\mu\nu}^{\lambda}$ are determined by the
Littlewood-Richardson rule (and are, in fact, interpretable as
Clebsch-Gordan coefficients for $SL\left(  n\right)  $). Because of the
triangular decomposition of the product of two Schur polynomials in terms of
other Schur polynomials, and the triangular decomposition of a Jack symmetric
polynomial (Schur polynomials in particular) in terms of the monomial
symmetric functions we can infer that
\begin{equation}
\left(  s_{\delta}\left( \mathbf{ x} \right)  \right)  ^{d}=
\sum_{\lambda\leq d\delta}
c_{\lambda
}^{\left(  \alpha\right)  }P_{\lambda}^{\left(  \alpha\right)  }\left(
x\right)  \label{10}%
\end{equation}
for suitable coefficients $c_{\lambda}^{\left(  \alpha\right)  }$ (Note that
\[
span\left\{  s_{\lambda}\mid \lambda\leq\delta\right\} 
 =span\left\{  m_{\lambda}\mid \lambda\leq\delta\right\}
=span\left\{  P_{\lambda}^{\left(  \alpha\right)  }\mid \lambda\leq\delta\right\}
\]
which follows immediately from condition (i) in the definition of the
$P_{\lambda}^{\left(  \alpha\right)  }$ and the linear independence of the
$P_{\lambda}^{\left(  \alpha\right)  }$ which, in turn, follows immediately from
orthogonality property (ii) in the definition of the $P_{\lambda}^{\left(
\alpha\right)  }$.)

We now quote a specialization of a result of Macdonald (\cite{Mac2}, pg. 386)
that is, in turn, derived by
applying a particular change of variables applied to a formula due to Gross
and Richards \cite{GR}, and Kadell \cite{Kad}.

\begin{Lemma} For $Re\left(  d\right)  >0$, $Re\left(  p\right)
>-1$,
\[
\int_{\mathcal{S}_{n}}P_{\lambda}^{\left(  \frac{2}{d}\right)  }\left(
\mathbf{x}\right)  \prod_{i=1}^{n}\left(  x_{i}\right)  ^{p}\left( \Delta\left(
\mathbf{x}\right) \right)  ^{d}d\mathbf{x}=\frac{1}{\Gamma\left(  a+1\right)  }v_{\lambda
}\left(  d\right)  \prod_{i=1}^{n}\Gamma\left(  \lambda_{i}+p+1+d\left(
n-i + 1\right)/2  \right)
\]
where%
\begin{align*}
a &  =\left|  \lambda\right| +n\left(p+1+d(n-1)/2 \right) \\
v_{\lambda}\left(  d\right)   &  =\prod_{1\leq i<j\leq n}\frac{\Gamma\left(
\lambda_{i}-\lambda_{j}+d\left(  j-i+1\right)/2  \right)  }{\Gamma\left(
\lambda_{i}-\lambda_{j}+d\left(  j-i\right)/2  \right)  } \quad .%
\end{align*}
\end{Lemma}

\bigskip

Now note that the second factor in the integrand on the right hand side of (\ref{Indp}) 
can be written as
\begin{eqnarray*}
\left(  \prod_{1\leq i<j\leq n}\left(  x_{i}^{2}-x_{j}^{2}\right)  \right)
^{d}&=&\left(  \prod_{1\leq i<j\leq n}\left(  x_{i}+x_{j}\right)  \right)
^{d}\left(  \prod_{1\leq i<j\leq n}\left(  x_{i}-x_{j}\right)  \right)
^{d}\\
&=&\left(  s_{\delta}\left( \mathbf{ x} \right)  \right)  ^{d}\left(  \Delta
\left( \mathbf{ x} \right)  \right)  ^{d} \quad .%
\end{eqnarray*}
And so
\[
I_{n,1,p}=\int_{\mathcal{S}_{n}}\left(  \prod_{i=1}^{n}x_{i}^{p}\right)
\left(  s_{\delta}\left( \mathbf{x} \right)  \right)  ^{d}\left(  \Delta\left(
\mathbf{x}\right)  \right)  ^{d}d\mathbf{x} \quad .
\]
We now employ the expansion (\ref{10}) with $\alpha=\frac{d}{2}$ to obtain
\begin{align*}
I_{n,d,p}  &  =\int_{\mathcal{S}_{n}}\left(  \prod_{i=1}^{n}x_{i}%
^{p}\right)  \left(  \sum{\lambda\leq d\delta}
c_{\lambda}^{\left(  \frac{2}%
{d}\right)  }P_{\lambda}^{\left(  \frac{2}{d}\right)  }\left( \mathbf{ x} \right)
\right)  \left(  \Delta\left( \mathbf{ x} \right)  \right)  ^{d}dx\\
&  =\sum_{\lambda\leq d\delta}
c_{\lambda}^{\left(  \frac{2}{d}\right)  }\int_{\mathcal{S}%
_{n}}\left(  \prod_{i=1}^{n}x_{i}^{p}\right)  P_{\lambda}^{\left(  \frac
{2}{d}\right)  }\left( \mathbf{ x} \right)  \left(  a_{\delta}\left( \mathbf{ x} \right)  \right)
^{d}dx\\
&  =\frac{1}{\Gamma\left(  a+1\right)  }\sum_{\lambda\leq d\delta}
c_{\lambda
}^{\left(  \frac{2}{d}\right)  }v_{\lambda}\left(  \frac{d}{2}\right)
\prod_{i=1}^{n}\Gamma\left(  \lambda_{i}+r+\frac{d}{2}\left(  n-i\right)
\right) \quad .
\end{align*}
Or, using the fact that $\left|\delta\right| = n(n-1)/2$. 
\begin{Lemma}
\begin{equation*}
I_{n,d,p}  = \frac{1}{\Gamma \left(n\left(p+d(n-1)\right)+1 \right)} \sum_{\lambda\leq d\delta}
C_{\lambda}\prod_{i=1}^{n}\Gamma\left(  \lambda_{i}%
+p+1+d\left(  n-i\right)/2  \right)  
\end{equation*}
where coefficients $C_{\lambda}$ depend only on $n$ and $d$.
\end{Lemma}

\begin{Lemma}
If $\lambda=\left[  \lambda_{1},\ldots,\lambda_{n}\right]  $ is a partition of
weight $d\left|  \delta\right|  $ such that $\lambda\leq d\delta$, then
\[
\lambda_{i}\geq\mu_{i}\equiv\left[  \frac{d\left(  n-i\right)  +1}{2}\right] \quad .
\]
\end{Lemma}

\emph{Proof. } Let $\mathcal{P}$ be the set of partitions $\lambda=\left[
\lambda_{1},\ldots,\lambda_{n}\right]  $ satisfying the criteria
\begin{align*}
\lambda_{1}  &  \geq\lambda_{2}\geq\cdots\geq\lambda_{n}\geq0\\
\sum_{i=1}^{n}\lambda_{i}  &  =\frac{d}{2}n\left(  n-1\right) \\
\sum_{j=1}^{i}\lambda_{j}  &  \leq\sum_{j=1}^{i}d\left(  n-j\right)
=ni-\frac{1}{2}i\left(  i+1\right)  =\frac{d}{2}i\left(  2n-i-1\right) \quad .
\end{align*}
Since the total weight of such a $\lambda$ is fixed, in order to minimize a
particular part $\lambda_{i}$ we need to arrange it so that the parts
$\lambda_{j}$ to the left of $\lambda_{i}$ are as large as possible while the
$\lambda_{j}$ to the right of $\lambda_{i}$ are also large as possible
(otherwise we could shift some of $\lambda_{i}$'s weight to the right). In
fact, if
\[
\mu_{i}\equiv\min_{\lambda\in\mathcal{P}}\lambda_{i}%
\]
we'll need
\[
\mu_{i}=\mu_{i+1}=\cdots=\mu_{n}%
\]
to make the $\mu_{j}$ to the right as large as possible and
\[
\mu_{j}=d\left(  n-j\right)  \qquad,\qquad j=1,\ldots,i-1
\]
for the $\mu_{j}$ \ to the left to be as large as possible. But for such a
minimizing configuration $\mu$, we must also have
\begin{align*}
\frac{d}{2}n\left(  n-1\right)   &  =\sum_{j=1}^{n}\mu_{j}=\sum_{j=1}%
^{i-1}d\left(  n-j\right)  +\left(  n-i+1\right)  \mu_{i}\\
&  =dn\left(  i-1\right)  -\frac{d}{2}i\left(  i-1\right)  +\left(
n-i+1\right)  \mu_{i} \quad . %
\end{align*}
Solving this for $\mu_{i}$ yields
\[
\mu_{i}=\frac{d}{2}\left(  n-i\right) \quad .
\]
However, if $d$ is odd then
\[
\mu_{i}=\frac{d}{2}\left(  n-i\right)
\]
will be an integer only when $n-i$ is even. In such a case, the first integer
larger than $\mu_{i}$ would be
\[
\frac{d}{2}\left(  n-i\right)  +\frac{1}{2} \quad .%
\]
Accounting for this circumstance, we can write
\[
\mu_{i}=\min_{\lambda\in\mathcal{S}}\lambda_{i}=\left[  \frac{d\left(
n-i\right)  +1}{2}\right]  \quad ,
\]
noting that the added $\frac{1}{2}$ is innocuous in the cases when $d$
is even or when $d$ is odd and $n-i$ is even. \qed

\begin{Theorem2}
For $d \in \mathbb{Z}_{\ge0}$,
\begin{equation}
I_{n,d,p}= \frac{1}{\Gamma\left(n\left(p+d(n-1)\right)+1\right)}\Phi\left(  p\right)  \prod_{i=1}^{n}\Gamma\left( p+1+d\left(
n-i\right)/2  \right) \label{3.2} 
\end{equation}
 $\Phi\left(  r\right)  $ is a polynomial of degree $\leq$ $\frac{d}%
{4}n\left(  n-1\right)  $ if $d$ is even, \ or of degree $\leq\frac{d}%
{4}n\left(  n-1\right)  -\frac{1}{2}\left[  \frac{n}{2}\right]  $ if $d$ is odd.
\end{Theorem2}

\emph{Proof.} Applying Lemmas 3.1 and 3.2 we have
\begin{align*}
I_{n,d,p}   &  =\frac{1}{\Gamma\left( n\left(p+d(n-1)\right)+1\right)  }%
\sum_{\lambda\leq d\delta}
C_{\lambda}\prod_{i=1}^{n}\Gamma\left(  \lambda_{i}-\mu
_{i}+\mu_{i}+p+1+\frac{d}{2}\left(  n-i\right)  \right) \\
&  =\frac{1}{\Gamma\left(n\left(p+d(n-1)\right)+1\right)  } \\
&\quad \times \sum_{\lambda\leq d\delta}
C_{\lambda}%
\prod_{i=1}^{n}\left(  \mu_{i}+p+1+\frac{d}{2}\left(  n-i\right)  \right)
_{\lambda_{i}-\mu_{i}}\Gamma\left(  \mu_{i}+p+1+\frac{d}{2}\left(  n-i\right)
\right) \\
&  =\frac{1}{\Gamma\left( n\left(p+d(n-1)\right)+1\right)  }\left(  \sum_{
\lambda\leq d\delta} C_{\lambda}%
\prod_{i=1}^{n}\left(  \mu_{i}+p+1+\frac{d}{2}\left(  n-i\right)  \right)
_{\lambda_{i}-\mu_{i}}\right) \\
& \qquad \times \left( \prod_{i=1}^{n}\Gamma\left(  p+1+d\left(
n-i\right)/2  \right) \right) \quad .
\end{align*}

This shows that $I_{n,d,p}$ has the form (\ref{3.2}) with 
\[
\Phi\left(p\right) =
 \sum_{
\lambda\leq d\delta} C_{\lambda}%
\prod_{i=1}^{n}\left(  \mu_{i}+p+1+\frac{d}{2}\left(  n-i\right)  \right)
_{\lambda_{i}-\mu_{i}} \quad .
\]

We'll now obtain the stated upper bounds for $\deg \Phi$ by determining the total degree of 
the products of Pochhammer symbols (as $\lambda$ ranges over the set of partitions $\le d\delta$).
\begin{align*}
\deg \Phi \le & \deg \prod_{i=1}^{n}\left(  \mu_{i}+p+1+\frac{d}{2}\left(  n-i\right)  
\right)_{\lambda_{i}-\mu_{i}} \\
&=\sum_{\iota=1}^{n}\left(  \lambda_{i}-\mu_{i}\right)  \\ 
&  =\sum_{i=1}^{n}\lambda_{i}-\sum_{i=1}^{n}\mu_{i}\\
&  =\frac{d}{2}n\left(  n-1\right)  -\sum_{i=1}^{n}\left[  \frac{d\left(
n-i\right)  +1}{2}\right] \quad .
\end{align*}
Now when $d$ is even, we have $\left[  \frac{d\left(  n-i\right)  +1}%
{2}\right]  =\frac{d}{2}\left(  n-i\right)  $ and so%
\begin{align*}
\sum_{\iota=1}^{n}\left(  \lambda_{i}-\mu_{i}\right)   &  =\frac{d}{2}n\left(
n-1\right)  -\sum_{i=1}^{n}\frac{d\left(  n-i\right)  }{2}=\frac{d}{2}n\left(
n-1\right)  -\frac{d}{2}n^{2}+\frac{d}{4}n\left(  n+1\right) \\
&  =\frac{d}{4}n\left(  n-1\right)  \qquad\left(  \text{ if }d\in
2\mathbb{Z}\right) \quad .
\end{align*}
When $d$ is odd, we have
\begin{align*}
\sum_{i=1}^{n}\left[  \frac{d\left(  n-i\right)  +1}{2}\right]   &
=\sum_{i=1}^{n}\frac{d}{2}\left(  n-i\right)  +\sum_{\substack{\left(
n-i\right)  \\odd}}\frac{1}{2}\\
&  =\frac{d}{4}n\left(  n-1\right)  +\sum_{\substack{\left(  n-i\right)
\\odd}}\frac{1}{2} \quad .%
\end{align*}
Note that if $n=2k$, then $\left(  n-i\right)  $ will be odd exactly $k$ times as
$i$ ranges from $1$ to $n$, and if $n=2k+1$, $n-i$ will be odd exactly $k$
times. We can thus write%
\[
\sum_{i=1}^{n}\left[  \frac{d\left(  n-i\right)  +1}{2}\right]  =\frac{d}%
{4}n\left(  n-1\right)  +\frac{1}{2}\left[  \frac{n}{2}\right]
\]
and so%
\[
\sum_{\iota=1}^{n}\left(  \lambda_{i}-\mu_{i}\right)  =\frac{d}{4}n\left(
n-1\right)  -\frac{1}{2}\left[  \frac{n}{2}\right]  \qquad,\qquad d\in\left[
1\right]  _{2} \quad . %
\]
\qed

\begin{Remark}
Explicit calculations reveal that for $d=1,2,4$ and $n=2,3,4,5$ that the stated upper
bounds on the degree of $\Phi\left(  r\right)  $ are, in fact, realized (and thus far without
exception). Thus,
$\frac{d}{4}n\left(  n-1\right)  $ is indeed the least upper bound for
general $d$ and $n$. 
We thus have a noncombinatorial proof of the following interesting and purely combinatorial
statement:%
\[
\sum_{\sigma_{1}\in\frak{S}_{n}}\cdots\sum_{\sigma_{d}\in\frak{S}_{n}%
}sgn\left(  \sigma_{1}\cdots\sigma_{d}\right)  \prod_{i=1}^{n}\Gamma\left(
r+2\sigma_{1}\left(  i\right)  +\cdots2\sigma_{d}\left(  i\right)  \right)
=\Phi\left(  r\right)  \prod_{i=1}^{n}\Gamma\left(  r+d\left(  n-i\right)
\right)
\]
where $\Phi\left(  r\right)  $ is a polynomial of degree $\le \frac{d}{4}n\left(
n-1\right)  $ (our point here being to highlight again the remarkable fact
that after removing the ``greatest common gamma factors factor'' from the
terms on the left hand side, exactly half of the leading terms cancel.)  
\end{Remark}

\section{{Recursive formulas}}

In this section, we demonstrate another evaluation tactic for Selberg-like integrals;
one aimed at exploiting the interplay between the integrand and the region of integration 
$\mathcal{S}_{n}$. We begin by considering integrals of the form
\begin{equation}
J_{n,\kappa}\left(f\right)  \equiv
\int_{\mathcal{S}_{n}}f \left(
\mathbf{x}\right) \left(\Delta^{(n)}(\mathbf{x}) \right)^{\kappa}  d\mathbf{x} \label{4-1}%
\end{equation}
where $f$ is some homogeneous symmetric function and $\mathcal{S}_{n}$ is the same
simplex as in (\ref{6}) and $\Delta^{(n)}(\mathbf{x})$ is the Vandermonde determinant in
$n$-variables.
\[
\Delta^{(n)}(\mathbf{x}) = \left(  \prod_{1\leq i<j\leq n}\left(  x_{i}-x_{j}\right)
\right) \quad . 
\]
Shortly, we shall specialize to the (still general) case when 
$f=\Phi_{\lambda}$, $\left\{ \Phi_{\lambda} \right\}$ being some standard
basis for homogeneous symmetric polynomials in $n$ variables indexed of partitions
$\lambda$.

\bigskip

\subsection{A Change of Variables Formula}
In what follows, it clarifies matters to indicate by $\mathbf{x}_{(n)}$ elements of $\mathbb{R}^{n}$,
and by $\mathbf{x}_{(i)}$ the vector in $\mathbb{R}^{i}$ obtained by dropping the last $n-i$
components of $\mathbf{x}_{(n)}$.

Consider the following change of variables
\begin{align}
t_{i}  &  =\frac{x_{i}-x_{n}}{1-nx_{n}}\qquad,\qquad x_{i}=\left(
1-t_{n}\right)  t_{i}+\frac{1}{n}t_{n}\label{12} \quad ,\\
t_{n}  &  =nx_{n}\quad\quad\ \qquad,\ \qquad x_{n}=\frac{1}{n}t_{n} \quad .\nonumber
\end{align}
It is easy to check that Jacobian of this transformation is
\[
\left(  \frac{\partial x}{\partial t}\right)  =\det\left[
\begin{array}
[c]{cccc}%
1-t_{n} & 0 & \cdots & \frac{1}{n}\\
0 & \ddots & \ddots & \vdots\\
0 & \cdots & 1-t_{n} & \frac{1}{n}\\
0 & \cdots & 0 & \frac{1}{n}%
\end{array}
\right]  =\frac{1}{n}\left(  1-t_{n}\right)  ^{n-1}%
\]
and the new region of integration is
\[
\left\{  \mathbf{t}\in\mathbb{R}^{n}\mid t_{1}\geq t_{2}\geq\cdots\geq
t_{n-1}\geq0\quad,\quad\sum_{i=1}^{n-1}t_{i}\leq1\quad,\quad\ 0\leq t_{n}%
\leq1\right\}  \approx\mathcal{S}_{n-1}\times\left[  0,1\right] \quad .
\]
We thus have

\begin{Lemma}%
\[
\int_{\mathcal{S}_{n}}f\left(  \mathbf{x}_{\left(  n\right)  }\right)
d\mathbf{x}_{\left(  n\right)  }=\int_{0}^{1}\frac{1}{n}\left(  1-t_{n}%
\right)  ^{n-1}\left(  \int_{\mathcal{S}_{n-1}}f\left(  \mathbf{x}_{\left(
n\right)  }\left(  \mathbf{t}_{\left(  n\right)  }\right)  \right)
d\mathbf{t}_{\left(  n-1\right)  }\right)  dt_{n}%
\]
where
\begin{align*}
\mathbf{x}_{\left(  n\right)  } &  \equiv\left(  x_{1},\ldots,x_{n}\right) \quad ,  \\
\mathbf{t}_{\left(  n\right)  } &  \equiv\left(  t_{1},\ldots,t_{n}\right) \quad , \\
x_{i}\left(  \mathbf{t}_{\left(  n\right)  }\right)   &  =\left(
1-t_{n}\right)  t_{i}+\frac{1}{n}t_{n}\quad ,\\
x_{n}\left(  \mathbf{t}_{\left(  n\right)  }\right)   &  =\frac{1}{n}t_{n}\quad .%
\end{align*}
\end{Lemma}

\begin{Remark}
If one iterates the change of variables formula, then after $n-1$ reductions,
one arrives at
\begin{eqnarray*}
\int_{\mathcal{S}_{n}}f\left(  \mathbf{x}_{\left(  n\right)  }\right)
d\mathbf{x}_{\left(  n\right)  } &=&\int_{0}^{1}\frac{1}{n}\left(
1-t_{n}\right)  ^{n-1}\left(  \int_{0}^{1}\frac{1}{n-1}\left(  1-t_{n-1}%
\right)  ^{n-2} \cdots  \right.  
 \\
&& \left. \left. \cdots \left(  \int_{\mathcal{S}_{1}}f\left(
\mathbf{x}\left(  \mathbf{t}\right)  \right)  dt_{1}\right) \cdots  \right)
dt_{n-1}\right)  dt_{n}\\
&=&\frac{1}{n!} \int_{0}^{1}\left(  1-t_{n}\right)  ^{n-1}\left(
\int_{0}^{1}\left(  1-t_{n-1}\right)  ^{n-2} \cdots \right. \\
&& \left. \left.  \cdots  \left(  \int_{0}%
^{1}f\left(  \mathbf{x}\left(  \mathbf{t}\right)  \right)  dt_{1}\right) \cdots
\right)  dt_{n-1}\right)  dt_{n} 
\end{eqnarray*}
where
\begin{align*}
x_{1}\left(  \mathbf{t}\right)   &  =\left(  1-t_{n}\right)  \left(
1-t_{n-1}\right)  \cdots\left(  1-t_{2}\right)  t_{1}+\frac{1}{2}\left(
1-t_{n}\right)  \cdots\left(  1-t_{3}\right)  t_{2}+\cdots+\frac{1}{n}t_{n} \quad ,\\
x_{2}\left(  \mathbf{t}\right)   &  =\frac{1}{2}\left(  1-t_{n}\right)
\cdots\left(  1-t_{3}\right)  t_{2}+\frac{1}{3}\left(  1-t_{n}\right)
\cdots\left(  1-t_{3}\right)  x_{3}+\cdots+\frac{1}{n}t_{n} \quad ,\\
&  \vdots\\
x_{n-1}\left(  \mathbf{t}\right)   &  =\frac{1}{n-1}\left(  1-t_{n}\right)
t_{n-1}+\frac{1}{n}t_{n} \quad ,\\
x_{n}\left(  \mathbf{t}\right)   &  =\frac{1}{n}t_{n} \quad .%
\end{align*}
If we adopt the convention that%
\[
\prod_{k=n+1}^{n}\left(  1-t_{k}\right)  \equiv1\qquad,
\]
the formula for $x_{i}\left(  \mathbf{t}\right)  $\ can be written a bit more
succinctly as
\[
x_{i}=\sum_{j=i}^{n}\frac{1}{j}\left(  \prod_{k=j+1}^{n}\left(  1-t_{k}%
\right)  \right)  t_{j} \quad .%
\]
It should be now clear that, via this cummulative change of variables,
every monomial $x_{1}^{m_{1}}\cdots x_{n}^{m_{n}}$ can be expressed as
a sum of terms of the form 
\[
ct_{1}^{a_{1}}\left(1-t_{1}\right)^{b_{1}} \cdots  t_{n}^{a_{n}}\left(1-t_{n}\right)^{b_{n}}
\]
and so the integral of a homogeneous symmetric polynomial $f\left(  \mathbf{x}_{(n)}\right)  
$ over $\mathcal{S}_{n}$ can be reduced to a sum of products of beta-integrals.
\end{Remark}
\begin{Remark}
The change of variables $\mathbf{x} \rightarrow \mathbf{x}
\left( \mathbf{t} \right)$ also facilitates the explicit computation of the integrals $I_{n,d,p}$.
For example, one finds in this way that
\begin{eqnarray*}
I_{2,1,p} &=&\frac{1}{\Gamma\left(2p+3 \right)}\left( \frac{2}{\sqrt{\pi }}%
\right) \Gamma \left( p\right) \Gamma \left( p+\frac{3}{2}\right) \quad , \\
I_{3,1,p} &=&\frac{  8p+15}{\Gamma \left( 2p+3\right) }{\Gamma \left( 3p+7\right) }\left( \frac{1}{2\sqrt{\pi }%
}\right) \Gamma \left( p\right) \Gamma \left( p+\frac{3}{2}\right) \Gamma
\left( p+2\right) \quad ,\\
I_{4,1,p} &=&\frac{\left(
16p^{2}+80p+93\right) }{\Gamma \left( 4p+13\right) }\left( \frac{2}{\pi }%
\right) \Gamma \left( p\right) \Gamma \left( p+\frac{3}{2}\right) \Gamma
\left( p+2\right) \Gamma \left( p+\frac{7}{2}\right) \quad , \\
I_{5,1,p} &=&\frac{\left(
99855+135232\,p+65456\,{p}^{2}+13568\,{p}^{3}+1024\,{p}^{4}\right)}{\Gamma \left( 5p+21\right) }\frac{3}{8\pi } \\
&&\times \Gamma \left( p\right) \Gamma \left( p+3/2\right) \Gamma \left(
p+2\right) \Gamma \left( p+7/2\right) \Gamma \left( p+4\right) \quad .
\end{eqnarray*}
Unfortunately, such explicit calculations have not given us any insight into the nature of the
polynomial factors $\Phi(p)$ (cf. Theorem 2).

\end{Remark}

\subsection{A recursive formula for $J_{n,\kappa}\left(
\Phi_{\lambda} \right)$ }

Now let $\left\{ \Phi_{\lambda} \right\}$ be some basis for the homogeneous
symmetric polynomials in $n$ variables and consider integrals of the form

\begin{equation}
J_{n,\kappa}\left( \Phi_{\lambda} \right)  \equiv
\int_{\mathcal{S}_{n}} \Phi_{\lambda} \left(
\mathbf{x}\right) \left(\Delta^{(n)}(\mathbf{x}) \right)^{\kappa}  d\mathbf{x} \quad .\label{2-11}% 
\end{equation}

We begin by noting that every symmetric basis function $\Phi_{\lambda}\left(
\mathbf{x}_{\left(  n\right)  }\right)  =\Phi_{\lambda}
\left(  x_{1},\ldots,x_{n}\right)  $ 
degree $\left|  \lambda\right|  $, will have a ``symmetric Taylor
expansion'' of the form
\begin{equation}
\Phi_{\lambda} \left(  \mathbf{x}_{(n)}
\right)  =\sum_{i=0}^{\left|  \lambda\right|  }\sum_{\left|  \mu\right|
=\left|  \lambda\right|  -i}c_{\lambda\mu}\Phi_{\mu}
\left(  \mathbf{x}_{\left(  n-1\right)  }\right)  x_{n}^{i}  \label{Phi}%
\end{equation}
for suitable coefficients $c_{\lambda\mu}$. We will refer to the
$c_{\lambda\mu}$ as the ``generalized Taylor coefficients'' for $\Phi
_{\lambda}$. Similarly, we can define ``generalized binomial
coefficients'' $\left(
\begin{array}
[c]{c}%
\lambda\\
\mu
\end{array}
\right)  $ for the $\Phi_{\lambda}$ by the formula%
\[
\Phi_{\mathbf{\lambda}}\left(  \mathbf{x}_{\left(
n\right)  }+t\mathbf{1}_{\left(  n\right)  }\right)  =\sum_{\mathbf{\mu}%
}\left(
\begin{array}
[c]{c}%
\lambda\\
\mu
\end{array}
\right)  \Phi_{\mathbf{\mu}}\left(  \mathbf{x}_{\left(
n\right)  }\right)  t^{\left|  \lambda\right|  -\left|  \mu\right|  } \quad . %
\]
In terms of these gadgets,\footnote{These ``gadgets'' are in fact interesting
mathematical objects in their own right (see, for example, \cite{VK} \S 2.5.1} 
we see that under the
substitutions%
\begin{align*}
x_{i}  &  =\left(  1-t_{n}\right)  t_{i}+\frac{1}{n}t_{n}=\left(
1-t_{n}\right)  \left(  t_{i}+\frac{t_{n}}{n\left(  1-t_{n}\right)  }\right)
\quad , \\
x_{n}  &  =\frac{1}{n}t_{n} \quad ,%
\end{align*}
the right hand side of (\ref{Phi}) becomes
\begin{eqnarray*}
&&\sum
_{i=0}^{\left|  \lambda\right|  }\sum_{\left|  \mu\right|  =\left|
\lambda\right|  -i}c_{\lambda\mu}\Phi_{\mu}\left(
\mathbf{x}_{\left(  n-1\right)  }\left(  \mathbf{t}_{\left(  n\right)
}\right)  \right)  x_{n}^{i}\left(  \mathbf{t}_{\left(  n\right)  }\right) \\
&&=\sum_{i=0}^{\left|  \lambda\right|  }\sum_{\left|  \mu\right|  =\left|
\lambda\right|  -i}c_{\lambda\mu}\Phi_{\mu}\left(
\left(  1-t_{n}\right)  \left(  \mathbf{t}_{\left(  n-1\right)  }+\left(
\frac{t_{n}}{n\left(  1-t_{n}\right)  }\right)  \mathbf{1}_{\left(
n-1\right)  }\right)  \right)  \left(  \frac{t_{n}}{n}\right)  ^{i}\\
&& =\sum_{i=0}^{\left|  \lambda\right|  }\sum_{\left|  \mu\right|  =\left|
\lambda\right|  -i}\left(  1-t_{n}\right)  ^{\left|  \mu\right|  }%
c_{\lambda\mu}\Phi_{\mu}\left(  \mathbf{t}_{\left(
n-1\right)  }+\left(  \frac{t_{n}}{n\left(  1-t_{n}\right)  }\right)
\mathbf{1}_{\left(  n-1\right)  }\right)  \left(  \frac{t_{n}}{n}\right)
^{i}\\
&&  =\sum_{i=0}^{\left|  \lambda\right|  }\sum_{\left|  \mu\right|  =\left|
\lambda\right|  -i}\sum_{\nu}\left(  1-t_{n}\right)  ^{\left|  \mu\right|
}c_{\lambda\mu}\left(
\begin{array}
[c]{c}%
\mu\\
\nu
\end{array}
\right)  \Phi_{\nu}\left(  \mathbf{t}_{\left(
n-1\right)  }\right)  \left(  \frac{t_{n}}{n\left(  1-t_{n}\right)  }\right)
^{\left|  \mu\right|  -\left|  \nu\right|  }\left(  \frac{t_{n}}{n}\right)
^{i}\\
&& =\sum_{i=0}^{\left|  \lambda\right|  }\sum_{\left|  \mu\right|  =\left|
\lambda\right|  -i}\sum_{\nu}\left(  n\right)  ^{\left|  \nu\right|
-n-i}\left(  t_{n}\right)  ^{\left|  \mu\right|  -\left|  \nu\right|
+i}\left(  1-t_{n}\right)  ^{\left|  \nu\right|  }c_{\lambda\mu}\left(
\begin{array}
[c]{c}%
\mu\\
\nu
\end{array}
\right)  \Phi_{\nu}\left(  \mathbf{t}_{\left(
n-1\right)  }\right)
\end{eqnarray*}
and
\begin{align*}
\left(\Delta\left( \mathbf{x}_{(n)}\right) \right)^\kappa &= \prod_{1\leq i<j\leq n}\left(  x_{i}-x_{j}\right)  ^{\kappa}  \\
& =\prod_{1\leq
i\leq n-1}\left(  x_{i}-x_{n}\right)  ^{\kappa}\prod_{1\leq i<j\leq
n-1}\left(  x_{i}-x_{j}\right)  ^{\kappa}\\
&  =\prod_{1\leq i\leq n-1}\left(  \left(  1-t_{n}\right)  t_{i}\right)
^{\kappa}\prod_{1\leq i<j\leq n-1}\left(  \left(  1-t_{n}\right)  \left(
t_{i}-t_{j}\right)  \right)  ^{\kappa}\\
&  =\left(  1-t_{n}\right)  ^{\kappa\left(  n-1+\frac{1}{2}\left(  n-1\right)
\left(  n-2\right)  \right)  }\left(  \prod_{1\leq i\leq n-1}t_{i}\right)
^{\kappa}\left(  \prod_{1\leq i<j\leq n-1}\left(  t_{i}-t_{j}\right)  \right)
^{\kappa}\\
&  =\left(  1-t_{n}\right)  ^{\kappa n\left(  n-1\right)  /2}\left(
\prod_{1\leq i\leq n-1}t_{i}\right)  ^{\kappa}\left(  \prod_{1\leq i<j\leq
n-1}\left(  t_{i}-t_{j}\right)  \right)  ^{\kappa} \quad .%
\end{align*}
Applying Lemma 4.1, we thus obtain
\begin{eqnarray*}
J_{n,\kappa}\left(  \Phi_{\lambda}\right)   &
=&\int_{\mathcal{S}_{n}}\Phi_{\lambda}\left(
\mathbf{x}_{\left(  n\right)  }\right)  \left(  \prod_{1\leq i<j\leq n}\left(
x_{i}-x_{j}\right)  \right)  ^{\kappa}d\mathbf{x}_{\left(  n\right)  }\\
&  =&\int_{0}^{1}\frac{1}{n}\left(  1-t_{n}\right)  ^{n-1} \\
&& \times \left(
\int_{\mathcal{S}_{n-1}}\Phi\left(  \mathbf{x}_{\left(  n\right)  }\left(
\mathbf{t}_{\left(  n\right)  }\right)  \right)  \prod_{1\leq i<j\leq
n}\left(  x_{i}\left(  \mathbf{t}_{\left(  n\right)  }\right)  -x_{j}\left(
\mathbf{t}_{\left(  n\right)  }\right)  \right)  ^{t}d\mathbf{t}_{\left(
n-1\right)  }\right)  dt_{n}\\
&  =&\sum_{i=0}^{\left|  \lambda\right|  }\sum_{\left|  \mu\right|  =\left|
\lambda\right|  -i}\sum_{\nu}\left(  n\right)  ^{\left|  \nu\right|
-n-i-1}c_{\lambda\mu}\left(
\begin{array}
[c]{c}%
\mu\\
\nu
\end{array}
\right) \\
&&  \times \left( \int_{0}^{1}\left(  t_{n}\right)  ^{\left|  \mu\right|  -\left|
\nu\right|  +i}\left(  1-t_{n}\right)  ^{n-1+\frac{\kappa}{2}n\left(
n-1\right)  +\left|  \nu\right|  } dt_{n} \right) \\
&& \times \int_{\mathcal{S}_{n-1}}\Phi_{\nu}
\left(  \mathbf{t}_{\left(  n-1\right)  }\right)  \prod_{1\leq
i<j\leq n-1}\left(  t_{i}-t_{j}\right)  ^{\kappa}d\mathbf{t}_{\left(
n-1\right)  } \quad .%
\end{eqnarray*}
Employing Euler's formula for the beta function%
\[
B\left(  r,s\right)  \equiv\frac{\Gamma\left(  r\right)  \Gamma\left(
s\right)  }{\Gamma\left(  r+s\right)  }=\int_{0}^{1}t^{r-1}\left(  1-t\right)
^{s-1}dt
\]
we thus obtain%
\begin{Theorem3}
\begin{eqnarray*}
J_{n,\kappa}\left(  \Phi_{\lambda}\right)
&=&\sum
_{i=0}^{\left|  \lambda\right|  }\sum_{\left|  \mu\right|  =\left|
\lambda\right|  -i}\sum_{\nu}\left(  n\right)  ^{\left|  \nu\right|
-n-i-1}c_{\lambda\mu}\left(
\begin{array}
[c]{c}%
\mu\\
\nu
\end{array}
\right) \\
&&\times  B\left(  \left|  \lambda\right|  +\left|  \nu\right|  +1,n-1+\frac
{\kappa}{2}n\left(  n-1\right)  +\left|  \nu\right|  \right)  J_{n-1,\kappa
}\left(  \Phi_{\nu}\right) \quad .
\end{eqnarray*}
\end{Theorem3}

\end{document}